\UseRawInputEncoding 
\documentclass[11pt]{amsart}
\usepackage{geometry}                % See geometry.pdf to learn the layout options. There are lots.
\usepackage[parfill]{parskip}    % Activate to begin paragraphs with an empty line rather than an indent
\usepackage{graphicx}
\usepackage{amssymb}
\usepackage{epstopdf}
\usepackage[all, 2cell]{xy}
\usepackage[section]{placeins}

\usepackage{hyperref}
\hypersetup{
    colorlinks=true,
    linkcolor=blue,
    filecolor=magenta,      
    urlcolor=cyan,
}

 \usepackage{breakurl}

\title{Voevodsky's Unachieved Project}

\author{Andrei Rodin}

\date{\today}

\begin{document}
\maketitle

\begin{abstract}
In a series of lectures given in 2003 soon after receiving the Fields Medal for his results in the Algebraic Geometry Vladimir Voevodsky (1966-2017) identifies two strategic goals for mathematics, which he plans to pursue in his further research. The first goal is to develop a ``computerised library of mathematical knowledge'', which supports an automated proof-verification. The second goal is to ``bridge pure and applied mathematics''. Voevodsky's research towards the first goal brought about the new Univalent foundations of mathematics. In view of the second goal Voevodsky in 2004 started to develop a mathematical theory of Population Dynamics, which involved the Categorical Probability theory. This latter project did not bring published results and was abandoned by Voevodsky in 2009 when he decided to focus his efforts on the Univalent foundations and closely related topics. In the present paper, which is based on Voevodsky's archival sources, I present Voevodsky's views of mathematics and its relationships with natural sciences, critically discuss these views, and suggest how Voevodsky's ideas and approaches in the applied mathematics can be further developed and pursued. A special attention is given to Voevodsky's original strategy to bridge the persisting gap between the pure and applied mathematics where computers and the computer-assisted mathematics have a major role. 
\end{abstract}

\section{Introduction}
Vladimir Voevodsky (1966-2017) made important contributions into two areas of mathematics, which until recently were thought of as barely related. One is the Algebraic Geometry, more specifically Motive theory. Voevodsky's contribution to this field in 2002 won him the Fields Medal. The other is Type theory and Foundations of Mathematics. In mid-2000s Voevodsky started a project of building new ``univalent'' foundations of mathematics, which support an automated proof-checking \cite{Grayson:2018}. This project boosted an ongoing research in Homotopy Type theory (discovered around 2006 independently by Voevodsky and Awodey\&Warren  \cite[Introduction]{UF:2013}), and put the automated proof-checking into a new theoretical perspective. Along with these great professional successes in the mid-2000s Voevodsky also experienced  what he later described  as ``the greatest scientific failure'' in his life (see Presentation 3 below). This unachieved  project aimed at a mathematical theory of Population Dynamics that could allow one, in particular, to reconstruct a history of human population on the basis of its present genetic constitution. Voevodsky actively worked on this project  during at least 5 years, from 2004 to 2009. In 2009 Voevodsky abandoned this applied project and until the end of his life and career focused all his efforts on the Univalent Foundations and related topics. 

Unpublished materials related to Voevodsky's project in the Population Dynamics  along with some other Voevodsky's unpublished materials are now publicly available via Vladimir's memorial page \url{http://www.math.ias.edu/Voevodsky/} maintained by Daniel Grayson in the Princeton Institute of Advanced Studies. References to these and some other related materials are given in this paper in the form [S$n$] where $n$ is the number of entry described below in the ``Sources'' section; these entries contain permalinks to slides, audios, videos and unpublished notes. Materials directly relevant to Voevodsky's work in the Population Dynamics are [S \ref{hist}, \ref{gen}, \ref{singletext}, \ref{miamitalks}]; for a related  research in the Categorical Probability theory see [S \ref{probatext}, \ref{proba}, \ref{catproba1}, \ref{catproba2}, \ref{paths}]. For early stages of Voevodsky's research in the Foundations of Mathematics, which are contemporary with his research in the Population Dynamics and Categorical Probability (before 2010), see [S \ref{form},  \ref{fmht}, \ref{hlambdalong}, \ref{hlambdashort}, \ref{typesys}]. The Univalent Foundation project in its stable mature form is presented in [S \ref{ufb}, \ref{ufp}, \ref{ufpro}, \ref{ufw},  \ref{omuf}]. Voevodsky's reflexions on his own work including his abandoned project in the Mathematical Biology are found in his interview given to his friend Roman Mikhailov in 2012 [S\ref{int}]. Fragments of this interview are reproduced below in English translation as Presentation 3 and Presentation 4.

In the remaining part of this paper I do two things. First, I present  Voevodsky's general views of mathematics and Voevodsky's strategy of strengthening links between the pure mathematics and the natural sciences. This part is fully based on the available Voevodsky's record, and involves no attempt on my part to develop or criticise his views and ideas. This material is relevant because the research project in the Population Dynamics just like the Univalent Foundations project has been strongly motivated by Voevodsky's general reflexions about mathematics and his vision of its desired future. Second, I provide a critical analysis of these views and ideas and suggest some ways in which they can be further developed. Mathematical details of Voevodsky's drafts related to his project in the Population Dynamics (including his work in the Categorical Probability) are not taken into account in the present paper.

\section{Place and role of mathematics according to Vladimir Voevodsky}
This Section consists of four Presentations based on Voevodsky's recorded materials. Presentation 1 is based on transparencies of his lecture given at the AMS-India meeting in Bangalore, Dec. 17-20, 2003, see [S\ref{outside}] (no audio- or video-recording of this lecture is available). Presentation 2 is based on transparencies and a written note related to two lectures given earlier in the same year in Wuhan University, China, see  [S\ref{wuhan1}, \ref{wuhan2}]). These two Presentations consist of my rendering or paraphrasing of Voevodsky's thoughts presented on the transparencies and in the written note; they include 
 a number of Voevodsky's original wordings. I take a liberty not to use quotations marks in such cases in order to make the text better readable. The diagrams used in these Presentations are made with Latex after Voevodsky's hand-written diagrams found in the transparencies.   

Presentations 3 and 4 are fragments of interview given by Vladimir Voevodsky to his friend Roman Mikhailov in July 2012 [S\ref{int}] translated into English by myself. In the Presentation 3 Voevodsky overviews his intellectual development and explains origins and motivations of his project in Population Dynamics and of the Univalent Foundations project. Presentation 4 includes some further Voevodsky's ideas concerning the pure mathematics, applied mathematics, and the natural sciences.

\subsection{Presentation 1: Mathematics and the Outside World [S\ref{outside}]} 

Mathematics is an integral \textemdash \ albeit very special \textemdash\ part of general problem-solving activity, which in its turn is a basic pre-scientific human condition. Various practical problems, which are conceptualised and approached with the common aka \emph{conventional} thinking are more effectively solved via the \emph{mathematical modelling} aka applied mathematics; the applied mathematics gives rise to the so-called pure mathematics. The pure mathematics 
develops via  (i) solving such \emph{external} problems, which come via the mathematical modelling, and also (ii) via formulating and solving its own \emph{internal} problems, that is, via proving and disproving various mathematical \emph{conjectures}. The interaction between conventional thinking, mathematical modelling and the pure mathematics, i.e., the flow of problems and their solutions, proceeds as shown at Fig. 1.:

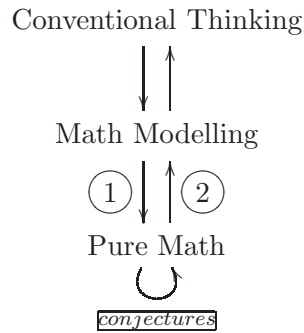
\begin{figure}[h!]
 $$\xymatrix{\textrm{Conventional Thinking}  \ar@<-5pt>[d] \\ \textrm{Math Modelling} \ar@<-5pt>[d]_*+[o][F]{1} \ar@<-5pt>[u]\\ \textrm{Pure Math} \ar@(dl,dr)[]_[F]{conjectures}\ar@<-5pt>[u]_*+[o][F]{2}}$$
\caption{Flow of Problems and Solutions}    
\end{figure}

Society supports mathematics, mainly, because of its capacity to solve problems arising in the applied mathematics using methods of pure mathematics (arrow 2 at the diagram) and also, eventually, for solving internal problems in the pure mathematics and for teaching old solutions to new generations.  

Over the last few decades the situation [as described above] was getting more and more out of balance. Arrows 1 and 2 shown at the diagram, which connect the pure and applied mathematics, were weakening. A weak incoming flow of external problems restrains the internal development of pure mathematics. A weak outgoing flow of useful solution restrains the support of mathematics provided by the Society. Breakdown of arrow 2 means eventually no salary for mathematicians. Breakdown of arrow 1  means eventually no new ideas in mathematics.    

How did we come to this poor situation? What can be done in order to improve it? What we need to do is to change the current pattern of using \underline{computers} in science. 

Presently computers enter into the above scheme of problem-solving as shown at Fig. 2:   

% \begin{figure}
% \includegraphics[scale = 0.1]{nfch1.pdf}

% \end{figure}
 
 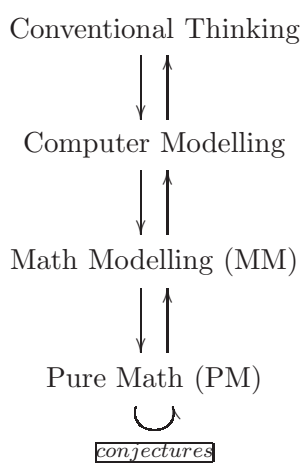
\begin{figure}[h!]
 $$\xymatrix{\textrm{Conventional Thinking}  \ar@<-5pt>[d] \\ \textrm{Computer Modelling} \ar@<-5pt>[d] \ar@<-5pt>[u]\\ \textrm{Math Modelling (MM)} \ar@<-5pt>[d] \ar@<-5pt>[u]\\ \textrm{Pure Math (PM)} \ar@(dl,dr)[]_[F]{conjectures}\ar@<-5pt>[u]}$$
\caption{The existing Flow Chart}    
\end{figure}
 
Here the flow of problems down to the ``mathematical modelling'' level is filtered through the ``computer modelling'' level. As a result  the ``mathematical modelling'' level, and as a consequence also the ``pure mathematics'' level, receive less problems than they used to receive before the rise of modern computer technologies. This particularly affects today's abstract mathematics. Problems, which pass through the filter, are formulated in the old-style language of variables and analytic functions, while the language of today's abstract mathematics is the Set theory. Thus at least a part of problems received at the ``pure mathematics'' level pass through a double-translation, which further weakens the incoming flow of external problems into the pure mathematics  (Fig. 3). 

% \begin{figure}
% \includegraphics[scale = 0.1]{dbltrans.pdf}
%\caption{Double Translation of Problems} 
 %\end{figure}
 
 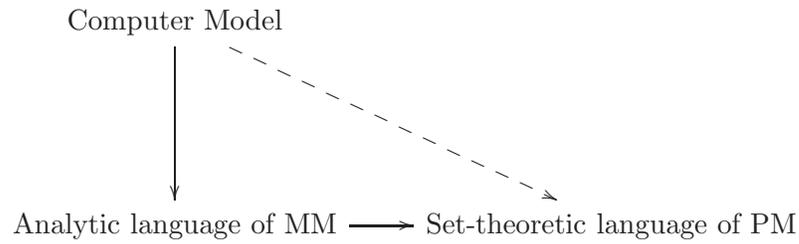
\begin{figure}[h!] 
$$\xymatrix{\textrm{Computer Model} \ar[dd] \ar@{-->}[ddr]\\ \\ \textrm{Analytic language of MM} \ar[r]& \textrm{Set-theoretic language of PM}}$$ 
\caption{Double Translation of Problems} 
 \end{figure}
 
The downstream flow of mathematical problems can be increased via a rearrangement of relationships between the computer modelling, on the one hand, and the pure mathematics, on the other hand, as shown at Fig. 4.:

%\begin{figure}
 %\includegraphics[scale = 0.1]{np.pdf}
%\caption{New Scheme of Relationships between the Computer Modelling and the Pure Mathematics} 
 %\end{figure} 

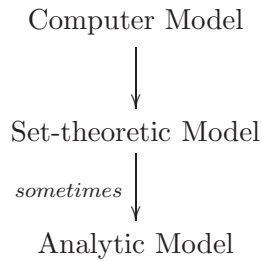
\begin{figure}[h!]
$$\xymatrix{ \textrm{Computer Model} \ar[d] \\ \textrm{Set-theoretic Model} \ar[d]_{sometimes} \\ \textrm{Analytic Model}}$$
\caption{New Scheme of Relationships between the Computer Modelling and the Pure Mathematics} 
 \end{figure} 

In order to implement this new scheme we need to reformulate fundamental and applied scientific theories in the language of today's abstract mathematics, viz., in the set-theoretic language. For this end we need to specify for each theory a notion of basic \emph{unit} and then consider sets of such units. Some examples are given in the Table 1.

\bigskip
\bigskip
\begin{table} [h] \centering
\begin{tabular}{|l|l|}
 \hline
 \textbf{Science} & \textbf{Unit} \\
  \hline	
  Population Biology and Demography & Individuals (individual organisms) \\
   \hline
   Financial Mathematics & Companies \\
   \hline
 Political Science & Voters \\
  \hline
  Particles Physics & Particles \\
    \hline
 Population Genetics & Genes \\
   \hline    \hline
 Future Theoretical Chemistry, which will be able \\ to account for individual molecules & Molecules \\
    \hline
 \end{tabular}
 
 \bigskip
 \caption{Scientific Disciplines and Their Ontological Units}
 \end{table}

The proposed rearrangement establishes a close connection between Science, Abstract Mathematics and Computing. It requires an essential reform in mathematical education.

The most important task for mathematicians is to produce examples that demonstrate the effectiveness of this approach.

\subsection{Presentation 2: What is most important for mathematics in the near future? Four Levels of Today's Mathematics.  [S\ref{wuhan1}, \ref{wuhan2}]}

There are two most urgent needs in today's mathematics:
\begin{enumerate}
\item To build a computerised library of mathematical knowledge, i.e., a computerised version of Bourbaki's \emph{Elements of Mathematics}; 
\item To bridge Pure and Applied Mathematics. 
\end{enumerate}  

[About (1).] We should gradually move from a hyperlinked mathematical text to a mathematical text verifiable with computer. 

 %[This is  all what Vladimir says in [S\ref{wuhan1}] about item (1). He also promises here to explain why Bourbaki's project, in his opinion, was not successful but doesn't elaborate. The rest of the note concerns item (2) and develops certain aspects of the topic treated in  [S\ref{outside}]. In Vladimir's later work the task (1) becomes central and leads him to the idea of Univalent Foundations, see   \underline{Presentation ??]} below. ]

[About (2).] Greatest advances of mathematics in the 20th century are in algebra, number theory and topology: they involve a combination of the visual intuition with the application of algebraic and symbolic methods. We discovered very fundamental classes of new  objects including categories, sheaves, cohomology, simplicial sets. They may turn out to be as important in science as algebraic groups. But presently we don't use them for solving problems outside the pure mathematics. 

One reason can be sociological. Only few people have a profound knowledge both of modern mathematics and of some other research  field where an application of modern mathematics can be possible. 

Another reason concerns the current scientific policies. In order to apply an abstract mathematical theory to a concrete practical problem one needs, first of all, to generalise this problem and abstract away the intuition associated with this problem. But the current funding policies favour rather fast solutions of concrete practical problems  such as, for example, designing ``the billion dollar drug''. 

In order to apply mathematics to a given problem outside mathematics one should begin with the opposite move.  Instead of trying to concentrate one's efforts on future applications of a mathematical theory to the real life, one should abstract oneself from the real life and look at the given problem as a formal game or puzzle. This is a reason why new mathematics too often strikes one, wrongly,  as detached from the real-world problems. 

So the only reasonable policy in mathematical research and in science in general is to support one's curiosity and one's sense of beauty in science.

 %Expression ``modern mathematics''  has two different senses. In the first sense it is a bulk of knowledge, i.e., a collection of theories, conjectures and theorems represented with texts. Most of these texts are purely mathematical while other texts provide connections between mathematics and other branches of knowledge. In the second sense the ``modern mathematics'' is a social institution, which consists of professional societies, mathematical departments and research institutes. Ultimately it consists of people.  

[About historical layers of today's mathematics.]
Today's mathematics comprises four different levels: 

\begin{enumerate}
\item \underline{Elementary Mathematics}: Pythagoras theorem, Quadratic Equations, etc.; it emerged more than 1000 years ago; 
\item \underline{``Higher'' Mathematics} : Integral and Differential Calculus, Differential Equations, Probability theory; it emerged in the 17th and 18th centuries; 
\item \underline{Modern Mathematics}: Modern Algebra (Galois theory, Group theory), Basic Topology, Logic (including G\"odel Incompleteness theorems) and Set theory; it emerged during the first half of  the 20th century; 
\item \underline{Synthetic Mathematics}: Representation theory, Algebraic Geometry, Homotopy theory (in particular the Motivic Homotopy theory), Differential Topology; it emerged during the second half of  the 20th century.
\end{enumerate}  

People who are not professional mathematicians usually know only levels (1) and (2), sometimes they have heard something about  level (3). The main reason is that 

\begin{itemize}
\item Elementary Mathematics is integrated into the everyday life;
\item ``Higher'' Mathematics is integrated into most sciences;
\item Modern Mathematics is integrated into \emph{some} sciences;
\item Synthetic Mathematics is very poorly integrated (if at all).
\end{itemize}

One of the most important challenges in mathematics today is to learn how to apply mathematics of higher levels in science and in the everyday life.

\subsection{Presentation 3: Origins and motivations of research in the Population Dynamics and in the Univalent Foundations [S\ref{int}, part 1].}

Since Fall 1997 I realised that my main contribution into the Motive theory and Motive Cohomology was already accomplished. Since then I was very consciously and actively looking for [$\dots$] a theme to work on after I accomplish my obligations related to Bloch-Kato conjecture\footnote{Earlier in the same interview Vladimir explains that during the last 5 years of his work in the Algebraic Geometry, which  stopped in 2010, he did not feel a real interest and continued his research only under the pressure of his self-imposed obligations.}.  [$\dots$] [C]onsidering tendencies of development of mathematics as a science I realised that we approach times when proving one more conjecture cannot change anything. [I realised] [t]hat mathematics is at the edge of crisis, more precisely, two crises. The first crisis concerns the gap between the ``pure'' and applied mathematics. It is clear that sooner or later there will arise the question of why the society should pay money to people, who occupy themselves with things having no practical application. The second crisis, which is less evident, concerns the fact that mathematics becomes very complex. As a consequence, once again, sooner or later mathematical papers will become too difficult for a detailed checking, and there will begin the process of accumulation of errors. Since mathematics is a very deep science in the sense that results of one particular paper usually depend on results of great many earlier papers, such an accumulation of errors is very dangerous for mathematics.

I decided to do something in order to prevent these crises. In the first case that meant to find an applied task, which would require for its solution methods of pure mathematics developed during the last years or at least during the last decades.

Since my childhood I was interested in natural sciences (physics, chemistry, biology) as well as in the theory of programming languages. Since 1997 I read a lot on these topics and took some university courses. In this way I significantly upgraded and deepened what I knew earlier. [$\dots$] [I] was looking for interesting open problems where I could apply today's mathematics. 
 
 Finally, I chose \textemdash\  as I now understand wrongly \textemdash\  the problem of reconstruction of history of populations [of living organisms] on the basis of their present genetic constitution. I worked on this problem for about two years
\footnote{Drafts [S\ref{hist}, \ref{gen}, \ref{probatext}]  are dated by 2004 (as the starting date of writing); in 2004 Voevodsky also gave a talk on the Categorical Probability theory motivated by his biological interests [S\ref{proba}]. These archival evidences show that in 2004 Voevodsky was already working on his project in the Population Dynamics, so this work continued at least 5 years in total.} 
 and finally realised in 2009 that everything that I invented was useless. This was the greatest scientific failure in my life. I great amount of work was invested into a project that failed completely. Some profit was gained anyway: I learned a lot about the Probability theory, which I didn't know well earlier, as well as about Demography and Demographic History.
 
 During the same time I was looking for an approach to the problem of accumulation of errors in works of pure mathematics. It was clear [to me] that the only solution was to create a language, with the help of which people can write mathematical proofs checkable with computer. Until 2005 this latter task seemed to me more difficult than the task in the historical genetics, on which I was was working. $\dots$ [However i]n 2005  I managed to formulate several ideas, which unexpectedly opened a possibility of new approach to one of the most important problems in the foundations of today's mathematics.

\subsection{Presentation 4: On the pure and applied mathematics; on the relevance of UF in the applied mathematics; on the data-driven science
 [S\ref{int}, part 2].}

Concerning the pure and applied mathematics, I have the following picture. Pure mathematics is working with models of high abstraction and low complexity (mathematicians like to call this low complexity elegance). Applied mathematics is working with more concrete models but on the higher complexity level (many equations, unknowns, etc). Interesting application of the modern pure mathematics are most likely in the area of high abstraction and high complexity. This area is practically inaccessible today, mostly due to the limitations of the human brain [$\dots$]. When we will learn how to use computers for working with abstract mathematical objects this problem will be no longer important and interesting applications of ideas of today's abstract mathematics will be found. 

That's why I think that my present work on computer languages
\footnote{Voevodsky refers here to his research in the Univalent Foundations and its implementation on computer.} that allow one to work with such objects, will be also helpful for application of ideas of today's pure mathematics in applied problems. 

\bigskip
[Then Voevodsky discusses difficulties of his project in the Population Dynamics related to insufficient empirical data and other non-mathematical matters.]

\bigskip
From a mathematical point of view the situation was also far from being ideal because nobody seriously studied such complex and temporally heterogeneous processes earlier. Finally, I came up with a new formalisation of Markov processes based on the notion of systems of paths. The paper turned out to be long and technical, and it still remains unfinished\footnote{By all evidence, Voevodsky refers here to [S\ref{probatext}].}
. I plan to return to this paper and finish it already with a convenient computer-based proof-assistant.  

\bigskip
[$\dots$]

\bigskip
Science should collect and comprehend a new knowledge. The collection part is very important. There is a view according to which all important observations are already done, the general world image is clear, so it remains only to arrange this knowledge and pack it into a compact and elegant theory. This view is wrong. It is not only wrong but also supports a very negative tendency to ignore everything that doesn't fit a ready-made theory or hypothesis. This is one of the most important problems of today's science.

\section{Discussion}

\subsection{General issues}
We can see from the Presentation 1 that Voevodsky's view on mathematics is broadly Aristotelian: he thinks of mathematics as an abstract form of reasoning, which is rooted in what he calls the \emph{conventional thinking}, i.e.,  a practice-oriented thinking that deals with multiple practical tasks\footnote{For a modern version of Aristotelian philosophy of mathematics see \cite{Franklin:2014}.}. \emph{Mathematical modelling}  mediates between the conventional thinking and the pure mathematics (Fig. 1 above). Mathematically-laden sciences are not distinguished in this scheme but the context makes it clear that the mathematical modelling also serves as an interface between science and pure mathematics.  Voevodsky grants for the pure mathematics a capacity to formulate and solve certain internal problems (the circular arrow at Fig. 1 ) but he claims that mathematics cannot successfully and sustainably develop without solving \emph{external} problems coming from the ``outside world''. He provides two independent reasons for this claim. 
 
One reason is sociological: mathematical research is supported by society primarily for its contribution to solving important practical problems. As Voevodsky makes it clear in the Presentation 2 his understanding of the social role of mathematics does not imply that all mathematical research should aim at solving concrete practical problems. In his view, a mathematical research can be more practically useful and effective when it is driven by one's curiosity and a sense of beauty rather than by practical interests.  %This reservation should not be interpreted in the sense that solving practical problems coming from the ``outside world'' plays no essential role in  mathematics.  

The second reason is epistemological: without a stable incoming flow of practical problems the pure mathematics soon runs out of new ideas. Thus we have here the following broad picture. Mathematics solves or help to solve practical problems by internalising them and abstracting them away from their local contexts. This triggers an internal dynamics in pure mathematics where certain problems are formulated and solved without further external motivations. However in a long term this internal mechanism is not sufficient for keeping a mathematical research going. So the incoming flow of external problem is just as vital for mathematics as the outgoing flow of practical solutions, which help mathematicians to secure a social status and necessary fundings for their research.

It is worth mentioning that Soviet philosophy of mathematics, which until 1990s remained relatively isolated from the Western mainstream along with the rest of Soviet philosophy, drew on the official Marxist doctrine of \emph{dialectical materialism} rather than the Analytic philosophical tradition that stressed logical and ontological foundations of mathematics.  As a result it gave more significance to applied and practical aspects of mathematics than to its logical and foundational aspects, see \cite{Ruzavin:1974}, \cite{Ruzavin:1983} by G.I. Ruzavin for standard examples. I don't know whether Vladimir Voevodsky ever read a philosophical literature on mathematics produced by Soviet philosophers but his overall picture of mathematics presented in  [S\ref{outside}, \ref{wuhan1}, \ref{wuhan2}] (Presentations 1, 2 above) is very similar to Ruzavin's: both authors describe mathematics in terms of its historical genesis from pre-scientific practical forms of problem-solving toward an internalised theoretical form of problem-solving. Even if Voevodsky's conviction that applications play a major role in the development of mathematics was not universally shared by all his  Russian colleagues, his view on this matter agree with an established local academic tradition represented by such prominent names as Andrey Kolmogorov, Vladimir Arnold and Israel Gelfand who consistently combined theoretical and applied mathematical research during all their professional careers. 

A systematic philosophical critique of Voevodsky's general view of mathematics is out of the scope of the present paper. In what follows I focus on specific details of this view, which concern functioning of channels between the pure mathematics and its applications. Recall that in Voevodsky's opinion during the last decades these channels have been seriously damaged  and now need an urgent repair.

\subsection{Pure and Applied Mathematics. The (in)effectiveness of con mathematics in the natural sciences.}
Unlike Eugene Wigner in 1960 \cite{Wigner:1960} Vladimir Voevodsky doesn't see the effectiveness of mathematics in the natural sciences as a miracle, and believes that mathematics can and should be applied in sciences more effectively  than it is presently applied. He finds it abnormal that the most recent mathematical achievements are very poorly integrated in the current science and technology, to leave alone the everyday practical life.
 
Even if such a view on today's mathematics is not unique\footnote{Israel Gelfand  famously remarked that  ``There is only one thing which is more unreasonable than the unreasonable effectiveness of mathematics in physics, and this is the unreasonable ineffectiveness of mathematics in biology'' \cite{Arnold:1998}.} it is also very far from being common.  According to an influential historical narrative the detachment of today's pure mathematics from science and technology is a consequence of a significant and positive conceptual change in mathematics and its foundations that occurred at the turn of the 20th century. Commenting in 1960 on Hilbert's 1899 \emph{Foundations of Geometry} \cite{Hilbert:1899}, \cite{Hilbert:1950} Hans Freudental remarks that 
with this Hilbert's achievement the ``bond with reality is cut. Geometry has become pure mathematics'' \cite[p. 618]{Freudental:1960}. The context makes it clear that by ``reality'' Freudental means here the reality of physical space and objects in this space but not a putative Platonic reality of self-standing mathematical objects  \cite{Balaguer:1998}. Freudental strongly praises this development and describes it retrospectively as a mainstream. In this context he quotes  Albert Einstein who says in his 1921 public lecture \emph{Geometry and Experience} (see \cite{Einstein:1922}) that ``the progress entailed by [the Hilbert-style \textemdash\ \emph{A.R.}] axiomatics consists of the sharp separation of the logical form and the realistic and  intuitive content''
\footnote{Quoted after \cite{Freudental:1960}; in the published version of this lecture \cite{Einstein:1922} the wording in this phrase is slightly different but its meaning is the same.} and that this modern conception of geometry plays a role in the Relativity Theory. Remarkably, Freudental does not mention that along with the axiomatic geometry  Einstein discusses in the same lecture the ``practical'' geometry, which in Einstein's view is at least as much important for science as the abstract axiomatic geometry. 
  
The idea according to which the pure mathematics is self-sustained and should be sharply separated from the applied mathematics and mathematically-laden sciences also motivated Bourbaki's long-term project of building new set-theoretic foundations of mathematics as well as  many related developments including the educational reform of school mathematics in the Bourbaki vein known in the US under the name of \emph{New Maths}. This educational reform started in the US in late 1950s (as a reaction on the launch of Soviet Sputnik), replicated in other countries,  but  by the early 1980s was already almost universally abandoned \cite{Phillips:2015}. Marshall Stone, a prominent mathematician who was a leading figure of the \emph{New Math} movement, expressed his understanding of the contemporary mathematics in his programmatic 1961 paper entitled  \emph{The Revolution in Mathematics} with the following strong claim:

\begin{quote}  
``While several important changes have taken place since 1900 in our conception of mathematics or in our points of view concerning it, the one which truly involves a revolution in ideas is the discovery that mathematics is entirely independent of the physical world.'' \cite[p.716]{Stone:1961}. 
\end{quote}  

In eyes of those people who share this view on the relationships between sciences and mathematics the fact that a large body of mathematical knowledge acquired after 1950 remains unapplied gives no reason for worries. It is readily explained and justified in terms of the ``new'' understanding of mathematics, which stresses its independence from the natural sciences and leaves its effectiveness in these sciences wholly unexplained and contingent. This provides a context in which one can wonder together with Eugene Wigner \cite{Wigner:1960} how mathematics can be possibly effective in physics and other natural sciences. A partial answer to (or at least a clarification of) this question is given by pointing to the historical fact that mathematical theories, which are effective in physics, have been built quite independently of their set-theoretic foundations and in many cases \emph{ab initio} involved physical motivations and intuitions.  

%(as evidenced by Eugene Wigner in 1960  \cite{Wigner:1960}).

 %\footnote{When Eugene Wigner in his 1960 paper \cite{Wigner:1960} calls the effectiveness of mathematics in the natural sciences ``unreasonable'' and describes it as a ``wonderful gift that we neither understand nor deserve'' he apparently assumes the same background understanding of the nature of mathematics.}. 

The idea of sharp separation between the logical form and the ``realistic content'' of scientific theories stressed by Freudental after Einstein back in 1960 was a cornerstone of the (very diverse and ramified) philosophical movement known as Logical Empiricism, which was very influential during the period of 1920-1950s but by the late 1960 was already seen by many important players as definitely dead \cite{Creath:2017}. This historical reference is sufficient for seeing that philosophical ideas and trends, however important and influential at certain point of history, should not be confused with scientific and mathematical achievements such as proofs of long-standing mathematical conjectures. The ``revolution in mathematics'' referred to by Stone in the above quote was indeed a revolution \emph{in ideas} about mathematics rather than a revolution in the mathematics itself. Unlike properly scientific revolutions such philosophical revolutions are not irreversible and not truly universal. The relationships between the pure mathematics and the natural sciences remains today, as ever, a matter of philosophical controversy, which involves many conflicting views and ideas.

\subsection{Set-theoretic foundations of mathematics and natural sciences}
 %The fact that the \emph{Elements of mathematics} by Bourbaki and all of Bourbaki-style set-based mathematics is associated with certain philosophical views on mathematics (which, of course, have many different varieties) does not make these views a proper part of this kind of mathematics. Nevertheless it has an important pragmatic consequence, which concerns the applicability of set-based mathematics in physics and other natural sciences. It goes without saying that

All mathematics  presently used in physics and other sciences can be  recast in the Bourbaki-style set-theoretic language and logically grounded on the set-theoretic foundations. Such a recasting  is shown at Voevodsky's \emph{double translation} scheme (see Fig. 2 and Fig. 3 in the Presentation 1 above), where the set-based \emph{pure mathematics} connects to \emph{computational models} and raw empirical data not directly but via the \emph{mathematical models} built with the standard toolkit of mathematical physics such as systems of partial differential equations and other  \emph{analytic} means (that usually date back to 19-18th century). These mathematical models function as an effective interface between the pure and applied mathematics and, at a larger scale, between mathematics and the ``outside world''. Theories of pure mathematics, which lay outside this interface, i.e., a large body of ``modern'' and more recent ``synthetic'' mathematics remain mostly unapplied. While in eyes of those mathematicians who share Hans  Freudental's and Marshall Stone's view such a situation is normal and even progressive and desirable, in eyes of Vladimir Voevodsky and other mathematicians who think of mathematics as an integral part of science, this is a harmful consequence of a strategic mistake that needs to be corrected. Among views of people of this second group Voevodsky's view on mathematics is distinguished by the fact the he, unlike the majority others  \cite{Arnold:1998}, \cite{Feynman:1965}, considers logical foundations of mathematics very seriously. In his 2003 Bangalore lecture (Presentation 1) Voevodsky identifies foundations of mathematics with the set-theoretic foundations.
Voevodsky changed this view after 2006 when he put forward the idea of alternative Univalent Foundations of mathematics [S \ref{fmht}, \ref{ufb}, \ref{ufp}, \ref{ufpro}].
 
Voevodsky's strategy to solve the double translation problem (as in 2003) is to built a shortcut from the set-based pure mathematics to  applied mathematics that would bypass the traditional ``analytic'' interface between the pure and applied mathematics. This proposal has an important computational aspect, which is discussed in \textbf{3.5} below. Here I focus on the ontological part of this proposal, which concerns a possibility to interpret abstract mathematical sets in naturalistic terms, see Table 1 in the concluding part of Presentation 1.        

I historical digression is here in order. The inventor of Set theory Georg Cantor back in 1884 conceived of possible applications of this theory in biology \cite{Ferreiros:2004}. The founder of modern axiomatic method David Hilbert conceived of possible applications of this method in natural sciences and included the problem of axiomatising physics into his famous list of 23 open mathematical problems that he announced in 1900 \cite{Hilbert:1902}. However the mainstream axiomatic Set theory since its very emergence in Zermelo's pioneering work \cite{Zermelo:1908} went in a wholly different direction, which did not involve such naturalistic considerations. This feature of abstract axiomatic Set theory was inherited by Bourbaki's set-theoretic foundations of mathematics. This is why Cantor's ideas concerning a possible role of Set theory in biology are commonly perceived today as a mere historical curiosity.  

Cantor's naturalistic speculations can be contrasted to Bernhard Riemann's speculations about a possible physical relevance of his notion of \emph{manifold} that bears today Riemann's name expressed in his 1854 Habilitation lecture \cite{Riemann:1868} and some other writings. Since the Riemannian geometry serves a mathematical foundation of the physical theory of General Relativity, which is nowadays universally accepted by the scientific community, it is justified to consider Riemann's geometrical studies and Einstein's work on GR as different stages of the same project and claim that Riemann's physical ideas pointed to the right direction. Given today's axiomatic Set theory, this is not the case of Cantor's naturalistic speculations. 
\footnote{The two cases are not independent: Cantor was influenced by Riemann's work and borrowed from him the term ``manifold'' (die Mannigfaltigkeit) that Cantor used for his set concept in his early works \cite{Ferreiros:1999}.}. 

The exceedingly abstract character of the standard mathematical set concept associated with ZF and akin axiomatic theories has been stressed and sharply criticised by William Lawvere in his seminal 1970 paper in the following words:

\begin{quote}  
``[A] 'set theory' $\dots$ should apply not only to \emph{abstract} sets divorced from time, space, ring of definition, etc., but also to more general sets, which do in fact develop along such parameters.''  \cite[p.329]{Lawvere:1970}\end{quote}  

\noindent Talking about ``general sets'' in the above quote Lawvere refers to the Topos theory but not to the Set theory in its standard axiomatic form. Attempts to re-establish the connection between the foundations of mathematics and the natural sciences, which in view of Hans Freudental, Marshall Stone and many other 20th century mathematicians has been lost irreversibly for a good reason, using Topos theory and Category theory will be briefly reviewed below in \textbf{3.4} 
  
Voevodsky's idea exposed in his Bangalore lecture [S\ref{outside}] (Table 1 in the Presentation 1 above) is not very original and it doesn't involve an attempt to revise the standard axiomatic foundations of Set theory. It amounts to a straightforward identification of mathematical sets with 
collections of material objects \textemdash\ sets of living organisms, sets of enterprises, sets of voters, sets of physical particles, sets of genes, and sets of molecules \textemdash\ and then using appropriate set-based mathematical structures in the corresponding scientific theories: general biology, economics, sociology, particle physics, genetics, chemistry and what not.  It is hardly possible to support the claim that this approach cannot work in principle but there are clear evidences that so far it has not been successful in spite of many attempts to apply it in the scientific practice.

The idea to use Bourbaki-style set-based formal representations of scientific theories was first proposed by Patrick Suppes and his collaborators back in the 1950s under the name of the ``semantic view of theories''. The  ``semantic view''  was proposed as a replacement for the ``syntactic view'', which these people ascribed to earlier enthusiasts of using new logical methods in science \cite{Halvorson:2016}. The earlier attempts to use logical methods in science included, on the one hand, some logically-based accounts of science offered by philosophers \cite{Nagel:1961} and, on the other hand, attempts to build workable scientific theories according logical recipes such as Joseph Henry Woodger's  axiomatic theory of biology \cite{Woodger:1937}. Proponents of the new ``semantic'' approach argued that a scientific theory cannot be identified with any particular axiomatic system but  should be identified instead with a class of models, which may, generally, satisfy different systems of axioms. In retrospect it is clear that Suppes and his followers continued earlier efforts to introduce logical methods into science but applied some more advanced logical techniques than their predecessors, which now included the Model theory and Formal Semantics. This formal approach had a significant impact in the Philosophy of Science of the last century \cite{Muller:2011}. However already at an earlier stage of this project it became obvious that even if the Bourbaki-style formal representation of scientific theories can be useful for certain philosophical and logical purposes (in particular, it may help to analyse logical relationships between different scientific theories and different branches of science), this way to represent theories is not appropriate for more common scientific purposes such as presentation of new theoretical results and teaching of university courses.  While the impact of Bourbaki's axiomatic style on the current mathematical practice remains significant albeit controversial, the impact of this axiomatic style in science is non-existent or negligible.

In Computer Science (CS) the idea to represent collections of material objects with abstract sets plays a role in  \emph{formal ontologies}, which provide a theoretical basis for computer-based Knowledge Representation technologies (KR) \cite{Brachman&Levesque:2004}. However KR in its existing forms does not reflect the theoretical structure of today mathematically-laden science such as theoretical physics and is not designed for this purpose. This is why it can hardly serve as an interface between science and set-based mathematics.

As I have already said, important fragments of set-based Bourbaki-style mathematics such as Group theory are successfully applied in physics and other sciences. However such mathematical theories are applied in sciences in different forms quite independently from  their set-theoretic representation and set-theoretic grounding. It is often said the modern axiomatic set concept is ``too abstract'' to be useful in science. I think that this claims needs qualifications, moreover that it is not clear how the degree of mathematical abstraction can be measured. Such mathematical concepts as that of natural and real number, of geometrical point, etc., are also abstract but nevertheless very useful and effectively applied in natural sciences. The identification of certain material objects, e.g., planets, with Euclidean geometrical points that ``have not part'' is a basic operation of the Classical Newtonian mechanics. By the same pattern one can     
%Such a straightforward identification of mathematical sets with collections of material objects may look, indeed, as a very suggestive, intuitive and natural idea. Easy objections such as ``the infinite sets do not exist in Nature'' are hardly convincing: points without parts, lengths without breadths,  infinitesimal quantities and infinite Cauchy series do not exist in Nature either but this doesn't prevent these mathematical concepts from being effectively applied in natural sciences via their identification with certain physical entities and processes\emph{modulo} some appropriate degree of abstraction. So emerge the concept of \emph{material point} in Classical Mechanics (that can be eventually be instantiated by a star or planet),  the concept of \emph{worldline} in Relativity Theory and similar physical concepts. Nevertheless a similar straightforward identification of an abstract mathematical set with a collection of physical particles or some other material system hardly plays a comparable role in science even if it has a strong intuitive appeal and can be often made by purpose. 
consider a set of molecules constituting a given living organism, and then assume that this set is a subject to the axioms of Zermelo-Fraenkel Set theory (ZF). However this way of thinking about a living organism or any other physical or biological system turns out to be scientifically sterile; it  doesn't make ZF relevant in science in anything like the same way in which Euclidean geometry is relevant in the Classical Mechanics, the Riemannian geometry is relevant in the General Relativity or operator algebras are relevant in the Quantum theory. 

 A plausible explanation of this situation can be in terms of historical epistemology can be given by pointing to the fact that unlike the Euclidean and Riemannian geometry, the axiomatic Set theory has not been developed in view of possible applications in natural sciences. The axioms of ZF are motivated by ideas that are metaphysical rather than physical (leaving now aside more specific mathematical and logical motivations). This concerns even such apparently unproblematic axioms as the \emph{axiom of pairing}. Informally, this axiom says that any given two things (objects, sets) $x, y$ can be always seen as one thing (pair) $\{x, y\}$. The intuitive appeal of this axiom is deceptive because our experience of collecting objects (or analysing a given complex object into its components) is always bound with certain spatial and temporal conditions, which the axiom of pairing wholly ignores. Motivations and alleged justifications of this mathematical axiom belong to the domain of speculative thought, which deliberately and systematically leaves aside all naturalistic, cognitive and empirical considerations  Whether or not this strategy is justified from an epistemological point of view, it obviously widens the gap between the pure mathematics and its applications\footnote{It worths mentioning that even when ZF is used a foundation of mathematics the axiom of pairing and other axioms of ZF are not applied in their ``official'' unrestricted form. The Bourbaki-style set-based mathematics involves an informal understanding based on implicit restrictions of type-theoretic character that certain sets and set-based structures are sound while some other, even if they are well-formed by the ZF standard, are not. As an example of such an unsound object think of set, elements of which are points of a given geometrical space and number $\pi$.}.

\subsection{Category theory as a mathematical foundation for natural sciences}
The above argument applies to the ``modern'' set-based Bourbaki-style mathematics but not to the more recent ``synthetic'' mathematics. In his Wuhan 2003 lectures [S\ref{wuhan1}, \ref{wuhan2}] Voevodsky does not make it quite clear what he means by the synthetic mathematics and how it differs from the ``modern'' mathematics (except providing some examples) but my guess is that he borrows term ``synthetic'' from William Lawvere and other people who use this for referring to category-theoretic axiomatic theories such as Synthetic Differential Geometry \cite{Kock:2010} and Synthetic Differential Topology \cite{Bunge:2018}. Lawvere and his followers reserve the term ``synthetic'' for their axiomatic category-theoretic approach in geometry, which they contrast to the more standard ``analytic'' approaches. One needs to bear in mind that in this specific context the term ``synthetic'' refers, by default, to theories that use Category theory as a foundation rather than to more more familiar Hilbert-style axiomatic theories that admit a set-theoretic semantics. When Voevodsky  in his 2003 lectures [S\ref{wuhan1}, \ref{wuhan2}] calls the contemporary Algebraic Geometry and Homotopy ``synthetic'' he, in my understanding, refers to this specific category-theoretic axiomatic approach applied in these mathematical disciplines \footnote{In his Paul Bernays Lectures delivered in Zurich in Fall 2014 [S\ref{eth}] Vladimir Voevodsky provides a more systematic and more detailed  account of historical development of mathematics and its foundations. This more detailed historical account agrees in its main features with that given in 2002 in the Wuhan lectures}.

Lawvere project of building new category-theoretic foundations of mathematics was strongly motivated by his wish to make the contemporary mathematics and its foundations more apt to applications in natural sciences. First attempts to use CT as a mathematical language for physical theories were made by Lawvere and a group of his followers in the early 1980s \cite{Lawvere:1982}; today there exists a significant number of accomplished works and ongoing research projects  that apply Category theory and Higher Category theory in various branches of mathematical physics \cite{Coecke&Paquette:2010}, \cite{Halvorson:2011}, \cite{Paugam:2014}
\footnote{For an updated survey and further references see also \url{https://ncatlab.org/nlab/show/higher+category+theory+and+physics}.}. In the theoretical biology similar attempts begin with Robert Rosen's pioneering 1958 paper \cite{Rosen:1958}. 

In these approaches the abstract mathematical concepts of category and functor are used for representing certain fundamental physical or biological concepts such as that of physical process and biological system. The resulting mathematically-laden scientific theories aim at representing fundamental principles and fundamental structures in their respective fields with new mathematical tools. Remarkably, Vladimir Voevodsky did not  enter into this vivid area of mathematical research but tried to  bridge today's mathematics with natural sciences in a different way.

\subsection{Applied Mathematics and Computers: Voevodsky's strategy}
As the above Presentation 4 makes it clear Vladimir Voevodsky was interested primarily in application of mathematics in an empirical 
data-driven scientific research, and opposed a tendency to reduce applications of new mathematics in science to fancy mathematical reformulations of earlier known theoretical results. This explains why Voevodsky did not invest his time and energy into mathematical physics, theoretical biology and akin theoretical disciplines.  In his 2003 Bangalore lecture [S\ref{outside}] (Presentation 1) Voevodsky stresses the role of electronic computers as an interface between between the pure mathematics, on the one hand,  and the ``outside world'', on the other hand (see Fig. 2 above). This role is twofold. First, the available computer memory allows for an extensive accumulation of empirical data\footnote{Since 2008 this phenomenon is commonly called the ``Big Data'' \cite{Editorial:2008}. An analysis of the impact of Big Data on scientific, social and political practices is found in \cite{Mainzer:2014}.}. Second, modern computers make it possible to proceed such large volumes of data. Nevertheless, as  Voevodsky argues in the Bangalore lecture,  the way in which computers are commonly used today in the data-driven sciences blocks the possibility to use the potential of contemporary mathematics in these sciences. This has a negative effect not only in these sciences but also in the pure mathematics because its poor interaction with the outside world deprives it from new ideas. 

Let me illustrate Voevodsky's argument with the example of Climate Research where electronic computers are systematically used both for storage and processing of empirical data and for computational modelling of climatic phenomena. The existing computer models of climate are typically based on mathematical models that involve mathematics of the pre-computer era: continuous functions, partial differential equations ( e.g., the Navier-Stokes equation) and akin analytic means. These computer models implement such traditional mathematical models via various numerical (computational) methods; in many cases this approach requires computational resources, which are at the limit of or exceed today's technological capacities  \cite{Palmer&Williams:2008}.

Since the mathematical foundation (say, the Navier-Stokes equation) of a computer model is fixed, the mathematical part of the task reduces to finding effective ways of solving the Navier-Stokes equation numerically. This can hardly produce any impact on the mainstream pure mathematics. Data and data structures that do not fit into this fixed mathematical framework are simply ignored and cannot motivate new developments in the pure mathematics.  Reciprocally, building computer models on the top of traditional mathematical analytic models leads to a waste of computational resources that cannot be fully compensated by inventing new sophisticated algorithms for finding numerical solutions without changing the theoretical foundation of these models. A revision of mathematical foundations of today's data-laden scientific theories such as the current theories of climate can help to use the available computational resources  more effectively.  As Voevodsky readily admits, until the effectiveness of this fundamental approach is demonstrated with concrete examples the above argument remains merely speculative. In order to provide such a working example Voevodsky focuses his efforts on a mathematical theory of Population Dynamics.

The above difficulty Voevodsky calls the problem of \emph{double translation} (see Fig. 3 above) for the following reason. In the Bangalore lecture he takes it for granted that the pure mathematics in its modern form resides on the set-theoretic foundations.  When scientific problems are translated into mathematical problems according the above pattern a mathematician needs, first, translate a given scientific problem (and relevant data structures) into an old-fashioned (analytic) mathematical language and, second, reformulate the problem in the modern mathematical language of Set-theory. Only after making these two preparatory steps one is in a position to  consider prospective applications of modern mathematics in the given area of scientific research. Voevodsky proposes to streamline this scheme by establishing a shortcut (shown with a dotted arrow at Fig. 3) from computer models to set-theoretics models. 

It should be born in mind that in 2003 when Vladimir Voevodsky delivered his Bangalore lecture he didn't conceive of any other possible foundation of the contemporary mathematics except the Set theory. Recall that in the same lecture he speculates about a possible role of Set theory as a mathematical foundation of natural sciences (see \textbf{3.3} above). In 2012 Voevodsky changes his ideas about foundations of mathematics and proposes alternative foundations,  namely, the Univalent Foundations (UF) (Presentations 3,4). UF was developed by Voevodsky in order to address the problem (``crisis'') of proof verification rather than the problem of application of modern mathematics (see Presentations 2 and 3). However UF, unlike the Set theory in its standard axiomatic form, also provides for the wanted shortcut from the computational mathematics to the foundations of today's mathematics. Talking in 2012 about his abandoned project in the Population Dynamics Voevodsky mentions his plan to accomplish this project (or at least its most developed part that involves the categorical probability theory [S\ref{probatext}) in the future using a proof-assistant (Presentation 4). Since in 2012 Voevodsky's research was fully focused on the UF, the timing makes it clear that mentioning a proof-assistant he points to UF here.

\subsection{Univalent Foundations and Applied Mathematics}
Since UF provides the wanted shortcut from computations to the contemporary ``synthetic'' mathematics (that includes Algebraic Geometry and Homotopy theory) and since such a shortcut, according to Voevodsky, helps to apply the contemporary mathematics in science and technology more effectively, it makes sense to consider UF as a possible mathematical foundation of scientific theories including biological theories. I can see at least two features of UF \textemdash\ one very general and the other more specific \textemdash\  which, in my view, make this mathematical theory an interesting candidate for this role. 

1) In low dimensions UF is supported with a spatio-temporal intuition that can be called homotopical intuition: the basic homotopical concepts of \emph{path} and \emph{path homotopy} are easily pictured and readily imagined \cite{Arkowitz:2011}. The significance of this feature does not reduce to its heuristic role. The homotopical intuition in the UF empowers an effective interface between the colloquial mathematical language used by working mathematicians, on the one hand, and formal proofs, on the other hand, thus making an automated verification of ``colloquial'' mathematical proofs possible. The same feature of UF may play a role in prospective applications of UF-based mathematics in natural sciences. As Ernest Cassirer put it in 1907 as a part of his critical assessment of Russell's mathematical logicism ``The principle according to which our concepts should be sourced in intuitions  means that they should be sourced in the mathematical physics and should prove effective in this field.'' \cite[p. 43]{Cassirer:1907} (my translation from German). The homotopical intuition allows one to think of Homotopy theory as a specific abstract form of  human spatio-temporal experience, which doesn't reduce to more familiar and better explored mathematical forms of experience associated with the standard arithmetical intuition, the traditional Euclidean geometrical intuition or the more modern Riemanian geometrical intuition. A dummy example that shows how homotopical concepts used in UF and the corresponding type-theoretic operations can be given a physical meaning is found in \cite{Rodin:2017}.

2) Unlike standard Hilbert-style axiomatic theories such as ZF, UF and UF-based mathematical theories are rule-based and may involve no axiom. This feature of UF facilitates its computational implementation
\footnote{The original version of UF uses rule-based Matin-L\"of Type theory (with dependent types, MLTT), a homotopical interpretation of  MLTT called Homotopy Type theory (HoTT), and an additional Axiom of Univalence (AU). However a more recent version of UF has a different formal carrier called Cubical Type theory (CTT) which is wholly rule-based and proves AU constructively as a theorem \cite{Cohen&Coquand&Huber&Mortberg:2016}}. Arguably, the same feature makes the UF-based constructive formal architecture more appropriate for scientific theories than the standard Hilbert-style axiomatic architecture. 

The homotopical semantics of type theories used in the UF distinguishes between propositional types (that represent logical propositions) and a hierarchy of higher non-propositional types, which represent various mathematical constructions that witness (are evidences for) their corresponding propositions. The presence of non-propositional types opens a possibility to interpret formal derivations in the UF-based formal theories as algorithms for building models with certain desired properties, in other words, for performing thought-experiments, which may justify and refute certain theoretical sentences. In the UF-based mathematics such thought-experiments translate into computational experiments, which may model certain physical processes. Thus the computable character of UF-based mathematical models suggests a possibility to use such models in the data-laden sciences such as the Climate Science where computational experiments play a major role. 

%Beware that the application of UF-based mathematics in science does not reduce to a mathematical trick, which could possibly help to develop faster algorithms for proceeding available empirical data and making more effective computations using the existing theoretical models. It rather amounts to a substantial revision and upgrade of the traditional apparatus of mathematical physics, which is still widely used today as a mathematical basis for building computational models in many areas of science.  

The possibility of such a direct application of today's abstract mathematical concepts to empirical data stored in a computer memory, which is suggested by the above analysis of Voevodsky's ideas, is independently explored in the Topological Data Analysis (TDA) \cite{Patania&Vaccarino&Petri:2017}  and its applications in the Neuroscience \cite{Expert:2019}, Biology \cite{Gameiro:2012}, \cite{Chan&Carlsson&Rabadan:2013} and some other scientific disciplines. This makes TDA an appropriate framework for a further development of Voevodsky's ideas presented in this paper. Among other contributions that these ideas could bring to TDA in its existing form is a theoretical computational basis and the proof-theoretic logical semantics (briefly outlined above) that UF inherits from MLTT.

\section{Conclusion}
The premature death in 2017 did not allow Vladimir Voevodsky to return to his project in applied mathematics, so it remained unachieved. However his general vision of today's mathematics and its possible future, which emphasises the interaction between the pure and applied mathematics, remains very inspiring and, as I hope, can motivate further research and lead to  important achievements. This concerns, in particular, Voevodsky's original strategy of bridging the growing gap between the pure and applied mathematics, which is focused on the empirical data-driven science rather than the theoretical science and aims at a fundamental renewal of the traditional apparatus of mathematical physics and mathematical biology.  In the previous sections of this paper I also tried to explain a possible relevance of the Univalent Foundations of mathematics developed by Voevodsky for a different purpose to his unachieved project in the applied mathematics. Vladimir Voevodsky's strategic vision of mathematics, including its past, its present and its projected future [S\ref{eth}]  is an important part of his intellectual heritage that needs a further study and further development.

\section*{Sources:}
\begin{enumerate}

\item \label{wuhan1} \textbf{What is most important for mathematics in the near future?} Wuhan University, November-December 2003, Lecture 1.\\
\underline{Permalink (notes):}  \\
\url{https://www.math.ias.edu/vladimir/sites/math.ias.edu.vladimir/files/2003_Wuhan_notes.pdf}

\item \label{wuhan2} \textbf{About mathematics and motivic homotopy theory } Wuhan University, November-December 2003, Lecture 2.\\
\underline{Permalink (transparencies):}  \\
\url{https://www.math.ias.edu/vladimir/sites/math.ias.edu.vladimir/files/2003_Wuhan_slides.pdf}

\item \label{outside} \textbf{Mathematics and the outside world}, Plenary lecture at AMS-India meeting Bangalore, December 17-20, 2003. \\ \underline{Permalink (transparencies):}   \\ \url{https://www.math.ias.edu/vladimir/sites/math.ias.edu.vladimir/files/Maths_and_the_outside_world.pdf}

\item \label{probatext} \textbf{A categorical approach to the probability theory}. February 3 - December 2, 2004. Unpublished unreleased manuscript
 \\ \underline{Permalink:}   \\ 
\url{https://www.math.ias.edu/Voevodsky//files/files-annotated/Dropbox/Unfinished_papers/Probability/Stochastic%20Categories/Stage2/stochastic.pdf}

\item \label{proba} \textbf{A categorical approach to the theory of probability} (in Russian), September 2, 2004. A talk in Steklov Institute (Moscow).\\ \underline{Permalink (announcement):}   \\ \url{http://www.mathnet.ru/php/seminars.phtml?presentid=1694}

\item \label{hist} \textbf{Historic inference from non-recombinant genetic data}, November 15, 2004 - March 6, 2006. 4 pages. An unfinished unreleased manuscript.
 \\ \underline{Permalink:}   \\ 
\url{http://www.math.ias.edu/Voevodsky/files/files-annotated/Dropbox/Unfinished_papers/Probability/Non%20Recombinant%20Genealogies/algorithms.pdf}

\item \label{gen} \textbf{What can the genetic data tell us about the history?},  April 14, 2004 - March 6, 2006. 3 pages. An unfinished unreleased manuscript
 \\ \underline{Permalink:}   \\ 
 \url{http://www.math.ias.edu/Voevodsky/files/files-annotated/Dropbox/Unfinished_papers/Probability/General%20Genealogies/histories.pdf}

\item \label{miamitalks} \textbf{Categories, Population Genetics and a Little of Quantum Physics}. A series of 3 lectures given in the University of Miami, Department of Mathematics in January 2005. Announcement and an abstract.
 \\ \underline{Permalink:}   \\ 
\url{https://www.math.miami.edu/highlights/events/lecture-series-by-v-voevodsky/}

\item \label{form} \textbf{On an approach to conveniently formalize mathematics}, June 4, 2005 - March 5, 2006. Unfinished and unpublished.  \\ \underline{Permalink:}   \\ 
\underline{Permalink:} \url{http://www.math.ias.edu/Voevodsky/files/files-annotated/Dropbox/Unfinished_papers/Dynamic_logic/Stage1/finitesets.pdf}

\item \label{fmht} \textbf{Foundations of Mathematics and Homotopy Theory}, March 22, 2006. IAS faculty lecture,  \\ \underline{Permalink (slides):}   \\ 
\url{https://www.math.ias.edu/vladimir/sites/math.ias.edu.vladimir/files/VV%20Slides.pdf}

\underline{Permalink (video):} \url{https://video.ias.edu/node/68}

\item \label{hlambdalong} \textbf{Notes on homotopy} $\lambda$\textbf{-calculus}. March 25, 2006. Unfinished and unpublished . \\ \underline{Permalink:}\\ 
\url{http://www.math.ias.edu/Voevodsky/files/files-annotated/Dropbox/Unfinished_papers/Dynamic_logic/Stage2/homotopy_lambda_calculus.pdf}

\item \label{hlambdashort} \textbf{A very short note on homotopy} $\lambda$\textbf{-calculus}. Preprint released in September 2006. \\ \underline{Permalink:}\\ \url{https://www.math.ias.edu/Voevodsky/files/files-annotated/Dropbox/For_Github/2006_09_A_very_short_note_on_homotopy_lambda/2006_09_Hlambda.pdf} 

\item \label{catproba1} \textbf{Categorical probability 1} (in Russian), November 20, 2008. A talk in the Steklov Mathematical Institute Seminar. \\ \underline{Permalink (announcement and video):}\\
\url{http://www.mathnet.ru/php/seminars.phtml?presentid=259}

\item \label{catproba2} \textbf{Categorical probability 2} (in Russian), November 25, 2008. A talk in the Steklov Mathematical Institute Seminar. \\ \underline{Permalink (announcement):}\\
\url{http://www.mathnet.ru/php/seminars.phtml?presentid=902}

\item \label{singletext} \textbf{Singletons}, May 11, 2009. 94 pages. An unfinished unreleased manuscript. 
\\
\underline{Permalink:}\\
\url{http://www.math.ias.edu/Voevodsky/files/files-annotated/Dropbox/Unfinished_papers/Probability/Singletons/Stage_2/singletons_current.pdf}

\item \label{paths} \textbf{Systems of paths} (in Russian), August 25, 2009. A talk in the  Steklov Mathematical Institute of RAS 
\\
\underline{Permalink (announcement and abstract):}\\
\url{http://www.mathnet.ru/php/seminars.phtml?presentid=752}

\item \label{typesys} \textbf{Notes on type systems}, September 8, 2009 - December 1, 2012. Unpublished
\\
\underline{Permalink:}\\
\url{http://www.math.ias.edu/Voevodsky/files/files-annotated/Dropbox/Unfinished_papers/Dynamic_logic/Stage_current/expressions_current-fixed.pdf}

\item \label{ufb} \textbf{Univalent Foundations}, Talk in Bonn, Sep. 8, 2010
\\
\underline{Permalink (Notes):}\\
\url{https://www.math.ias.edu/vladimir/sites/math.ias.edu.vladimir/files/Bonn_talk.pdf}
\\
\underline{Permalink (Coq examples):}\\
\url{https://www.math.ias.edu/vladimir/sites/math.ias.edu.vladimir/files/Bonn_talk_coq.pdf}

\item \label{ufp} \textbf{Univalent Foundations}, Talk at IAS, Dec. 10, 2010
\\
\underline{Permalink (video):}\\
\url{https://video.ias.edu/univalent/voevodsky}

\item \label{ufpro} \textbf{Univalent Foundations Project}. A modified version of an NSF grant application, October 2010
\\
\underline{Permalink:}\\
\url{https://www.math.ias.edu/vladimir/sites/math.ias.edu.vladimir/files/univalent_foundations_project.pdf}

\item \label{ufw} \textbf{Univalent Foundations}, Plenary lecture at WoLLIC, May. 18, 2011
\\
\underline{Permalink (slides):}\\
\url{https://www.math.ias.edu/vladimir/sites/math.ias.edu.vladimir/files/2011_WoLLIC.pdf}
\\
\underline{Permalink (Coq examples):}\\
\url{https://www.math.ias.edu/vladimir/sites/math.ias.edu.vladimir/files/demo_1.v}

\item \label{int} \textbf{Interview with Roman Mikhailov} (in Russian), July 2012
\\
\underline{Permalink (Part 1):}\\
\url{https://baaltii1.livejournal.com/198675.html}
\\
\underline{Permalink (Part 2):}\\
\url{https://baaltii1.livejournal.com/200269.html}

\item \label{omuf} \textbf{The Origins and Motivations of Univalent Foundations}, Article in the Institute Letter (IAS Princeton) written after the Lecture Mar. 26, 2014; appeared in Summer 2014
\\
\underline{Permalink:}\\
\url{https://www.ias.edu/ideas/2014/voevodsky-origins}

\item \label{eth} \textbf{Foundations of mathematics \textemdash\ their past, present and future.} The 2014 Paul Bernays Lectures at ETH Zurich. Sep. 9-10, 2014. 

\bigskip
\underline{Lecture 1}: \textbf{To the history of the conception}\\
\underline{Permalink (slides):}\\
\url{https://www.math.ias.edu/vladimir/sites/math.ias.edu.vladimir/files/2014_09_Bernays_1%20presentation.pdf} \\
\underline{Permalink (video):}\\
\url{https://www.math.ias.edu/vladimir/sites/math.ias.edu.vladimir/files/20140909_HGF5_Paul_Bernays_Lectures_Voevodsky_01-dm.mp4} 

\bigskip
\underline{Lecture 2}: \textbf{The story of Set theory (so far)}\\
\underline{Permalink (slides):}\\
\url{https://www.math.ias.edu/vladimir/sites/math.ias.edu.vladimir/files/2014_09_Bernays_2%20presentation.pdf}
 \\
\underline{Permalink (video):}\\
\url{https://www.math.ias.edu/vladimir/sites/math.ias.edu.vladimir/files/20140910_HGF5_Paul_Bernays_Lectures_Voevodsky_02-dm.mp4} 

\bigskip
\underline{Lecture 3}: \textbf{Univalent Foundations} \\
\underline{Permalink (slides):}\\
\url{https://www.math.ias.edu/vladimir/sites/math.ias.edu.vladimir/files/2014_09_Bernays_3%20presentation.pdf}
 \\
\underline{Permalink (video):}\\
\url{https://www.math.ias.edu/vladimir/sites/math.ias.edu.vladimir/files/20140910_HGF5_Paul_Bernays_Lectures_Voevodsky_03-dm.mp4}

\end{enumerate}

\bibliographystyle{plain} 
\bibliography{voevo1} 

\end{document}